\documentclass[3p,10pt]{elsarticle}
\usepackage{amsmath,amssymb,bm,graphicx,subfigure}
\usepackage{times}
\usepackage{diagbox}
\usepackage[colorlinks=true]{hyperref}
\usepackage{algorithmic, algorithm}
\usepackage{textcomp}
\usepackage{xcolor}
\usepackage{mathtools}
\usepackage{microtype}
\usepackage{multirow}
\newcommand{\tabref}[1]{Table~\ref{#1}} 
\newcommand{\figref}[1]{Fig.~\ref{#1}} 
\renewcommand{\eqref}[1]{Eq.~(\ref{#1})}
\newcommand{\secref}[1]{Section~\ref{#1}}
\newcommand\filledcirc{{\color{gray}\bullet}\mathllap{\color{gray}\circ}}
\DeclareMathOperator*{\minimize}{minimize}
\biboptions{sort&compress}
\begin{document}
\begin{frontmatter}
\title{Data-driven Experimental Modal Analysis by Dynamic Mode Decomposition}
\tnotetext[t1]{
This manuscript is the accepted version of the article:
A. Saito and T. Kuno,
``Data-driven Experimental Modal Analysis by Dynamic Mode Decomposition,''
Journal of Sound and Vibration, 481, 115434 (2020).
\url{https://doi.org/10.1016/j.jsv.2020.115434}
\copyright\ 2020 The Author(s). Published by Elsevier Ltd. This manuscript version is made available
under the CC-BY-NC-ND 4.0 license
(\url{https://creativecommons.org/licenses/by-nc-nd/4.0/}).
}
\author[meiji]{Akira~Saito\corref{cor1}}
\ead{asaito@meiji.ac.jp}
\author[meiji]{Tomohiro~Kuno}
\address[meiji]{Meiji University, Kawasaki, Kanagawa 214-8571, Japan}
\cortext[cor1]{Corresponding author. Tel:+81 44 934 7370.}
\begin{abstract}
This paper discusses the application of Dynamic Mode Decomposition (DMD) to the extraction of modal properties of linear mechanical systems, i.e., experimental modal analysis (EMA). First, theoretical background of the DMD is briefly reviewed and its relevance to the Ibrahim time-domain method is discussed. Second, DMD is applied to a single DOF system and multi-DOF discrete system to discuss the applicability and interpretation of the DMD as a method of EMA. Furthermore, the effects of measurement errors on the results of DMD are discussed.  It is shown that with relatively small measurement errors, DMD can capture modal parameters accurately. However, with relatively large measurement errors, DMD fails to capture modal parameters. 
Finally, DMD is applied to experimentally-obtained displacement field of a cantilevered beam, and its modal parameters are extracted. It is shown that the modal parameters extracted by DMD are as accurate as the ones obtained by the existing modal parameter extraction method. 
\end{abstract}
\end{frontmatter}
\section{Introduction}\label{section:introduction}
There have been many attempts to date to reveal modal properties of engineered structures that are subject to dynamic loadings. The modal properties of the structures are useful because their dynamics can be well described by a linear combination of the modes, when the dynamics of the structures resides in the linear dynamics regime. 
That enables ones to design, control, and optimize the dynamic properties of the structures accurately and appropriately. 

The modal properties are typically extracted from experimentally obtained frequency response functions~(FRFs) or time responses. 
This process is referred to as experimental modal analysis~(EMA) and is one of the key technologies in structural dynamics. 
The modal parameter extraction methods can roughly be classified into two categories: the ones based on frequency-domain signals, and the ones based on time-domain signals. 
Nowadays, it is more common to apply frequency-domain methods than the time-domain methods. 
In particular, a class of frequency domain methods that are based on least squares complex frequency (LSCF) domain method~\cite{VanDerAuweraerEtAl2001,PeetersEtAl2004} are widely used for its excellent capability to extract modal properties of various structures. 
The application of time-domain methods is limited to specialized applications especially in cases of structures with very low natural frequencies for which the time required to acquire several cycles worth of data can become a practical problem~\cite{Ewins2009}. 
%
%

In the area of fluid mechanics, for the purpose of extracting coherent structures from flow field, a dimensionality-reduction method called proper orthogonal decomposition~(POD) was proposed~\cite{BerkoozEtAl1993}. 
Later, the POD has been applied to nonlinear mechanical systems and widely used in structural dynamics area~\cite{FeenyKappagantu1998,KerschenEtAl2005}. 
Even though the importance of the modes extracted by the POD, which are referred to as POD modes, can be determined by the magnitude of the values of the corresponding eigenvalues, which are sometimes referred to as POD values, their {\it temporal} information, i.e, frequency and damping ratio, is completely lost. 
Recently, to circumvent this drawback of the POD, {\it dynamic mode decomposition}~(DMD) has been proposed by Schmid~\cite{Schmid2010} as an extension to POD which well suits a data-driven scenario where no governing equation is known. 
Unlike POD, the DMD can extract not only the modes, but also their temporal properties, i.e., frequencies and decay rates. 
Due to this property, DMD has been actively developed and applied to various applications involving nonlinear processes either numerically or experimentally, such as 
%
acoustic mode identification in a three dimensional chamber~\cite{JourdainEtAl2013}, a swirling flow problem~\cite{BistrianEtAl2017}, 
combustion~\cite{RichecoeurEtAl2012}, and 
pressure sensitive paint data~\cite{AliEtAl2016}, and many others. 

Since DMD can extract spatially coherent structures and their temporal properties from data, 
if the data is obtained from vibration tests of structures, these quantities correspond to modal parameters, i.e., mode shapes, natural frequencies, and damping ratios of the structures. 
To date, however, there have been few attempts to apply DMD for the structural dynamics applications, which may be attributed to the fact that the vibration modes are readily available in structural dynamics for representing coherent structures in the measured vibration data, both from numerically or experimentally obtained data. 
{\color {black}{From the viewpoint of EMA, since DMD is based on time-domain data, the method is classified into the time-domain method. Since the method has been actively developed for the extraction of modes for very large dimensional systems such as flow field, it has the potential to be applied to data-driven structural dynamics problems with a massive number of measurement points. With the increasing demand for capturing the distribution of accurate structural response over structures with many sensors, such capability of DMD would be beneficial.}}
In this paper, DMD is applied to extract modal properties of mechanical systems and discuss its capability and potential as a method of EMA. 

The remainder of the paper is organized as follows. In \secref{section:theory}, mathematical background of the DMD is briefly reviewed. 
In \secref{section:numerical_example}, numerical examples are provided by using discrete mass-spring-damper systems. 
In \secref{section:experimental_example}, DMD is applied to experimentally-obtained displacement field of vibrating cantilevered beam subjected to impulse loading, and the results are compared with those obtained by a conventional EMA method. 
In \secref{section:conclusion}, concluding remarks are provided. 
\section{Theory}\label{section:theory}
This sections briefly reviews the mathematical background of DMD from the viewpoint of structural dynamics. 
In particular, a focus is placed upon the relevance of the DMD to Ibrahim time domain (ITD) method~(see, e.g., \cite{Ewins2009}), which is a classical time-domain modal extraction method.  
\subsection{Dynamic Mode Decomposition}
Mathematical definition of the DMD is stated as follows. 
Let ${\bf x}(t)\in\mathbb{R}^m$ denote a vector of generalized coordinates, or a state vector of an $m$-dimensional dynamical system with $t$ denoting time. Consider a {\it snapshot matrix} ${\bf X}$ containing snapshots of the generalized coordinates ${\bf x}_j\triangleq{\bf x}(t_j)$ for $j=1,\dots,n$ as its column vectors, i.e., ${\bf X}=[{\bf x}_1,\dots,{\bf x}_n]$ for $t_{j+1}=t_{j}+\Delta t$. Even though DMD does not require $\Delta t$ be constant, we shall assume that it is constant throughout the paper. 
Furthermore, define a time-shifted snapshot matrix ${\bf Y}=[{\bf x}_2,\dots,{\bf x}_{n+1}]$. 
We then assume that there exists a linear transformation between ${\bf x}_j$ and ${\bf x}_{j+1}$, which is represented by a matrix ${\bf A}$, i.e., 
\begin{equation}
	{\bf x}_{j+1}={\bf A}{\bf x}_{j}\label{eq:th:-1}, 
\end{equation}
which leads to 
\begin{equation}
{\color {black}{
	{\bf Y} = {\bf A}{\bf X}, 
}}
\end{equation}
where 
\begin{equation}
{\bf A}={\bf Y}{\bf X}^{\dagger}\label{eq:th:-1+1}, 
\end{equation}
and ${\bf X}^{\dagger}$ denotes the Moore-Penrose pseudo inverse matrix of ${\bf X}$. 
The {\it DMD modes and eigenvalues} are defined as the eigenvectors and eigenvalues of ${\bf A}$. 
If singular value decomposition~(SVD) is applied to {\bf X}, then ${\bf X}={\bf U}\bm{\Sigma}{\bf V}^*$ where ${\bf U}\in\mathbb{C}^{m\times r}$, $\bm{\Sigma}\in\mathbb{R}^{r\times r}$, and ${\bf V}\in\mathbb{C}^{n\times r}$, $r$ is the rank of ${\bf X}$, and $^*$ denotes Hermitian transpose of a matrix. 
Note that by definition, ${\bf U}$ contains the {\it POD modes}. 
The pseudo inverse can then be computed as follows: 
\begin{equation}
{\bf X}^{\dagger}=
{\color{black}{{\bf V}\bm{\Sigma}^{-1}{\bf U}^*. 
}}
\end{equation}
The DMD eigenvalue are then obtained by computing the eigenvalues of the projected matrix ${\bf A}$ onto the POD modes~\cite{TuEtAl2014}. Namely, defining the projected matrix as 
\begin{equation}
\tilde{\bf A}\triangleq{\bf U}^*{\bf A}{\bf U}=
{\color{black}{
{\bf U}^*{\bf Y}{\bf V}\bm{\Sigma}^{-1}{\bf U}^*{\bf U}=
}}
{\bf U}^*{\bf Y}{\bf V}\bm{\Sigma}^{-1}, 
\end{equation}
the DMD eigenvalues are the eigenvalues of $\tilde{\bf A}$, i.e., 
\begin{equation}
\tilde{\bf A}{\bf w}=\mu{\bf w}, 
\end{equation}
where $\mu$ here is the DMD eigenvalue. 
The DMD mode corresponding to $\mu$ is then defined by 
\begin{equation}
\bm{\varphi}={\bf U}{\bf w}. \label{eq:th:32}
\end{equation}

{\color {black}{
Let us now discuss how damping ratios and natural frequencies can be obtained from the DMD eigenvalues and DMD modes from the viewpoint of linear systems theory. 
Denoting $i$-th DMD eigenvalue and mode respectively as $\mu_i$ and $\bm{\varphi}_i$, the following relationship holds:
\begin{equation}
\bm{\Lambda}=\bm{\Phi}^{*}{\bf A}\bm{\Phi}, 
\end{equation}
where $\bm{\Phi}=[\bm{\varphi}_1,\dots,\bm{\varphi}_m]$ and $\bm{\Lambda}={\rm diag}_{i=1}^m(\mu_i)$.  Applying the modal coordinate transformation ${\bf x}_{j}=\bm{\Phi}{\bf q}_j$ to \eqref{eq:th:-1} where ${\bf q}_j$ denotes the vector of modal coordinates at $t_j$, we obtain
\begin{equation}
{\bf q}_{j+1}=\bm{\Phi}^*{\bf A}\bm{\Phi}{\bf q}_{j}=\bm{\Lambda}{\bf q}_j. \label{eq:rev:1}
\end{equation}
Noting that \eqref{eq:rev:1} can be repeated from $j=k$ down to $j=0$, we obtain the following relationship, 
\begin{equation}
{\bf q}_k=\bm{\Lambda}^k{\bf q}_0, 
\end{equation}
where ${\bf q}_0=\bm{\Phi}^{*}{\bf x}_0$. Therefore, 
\begin{equation}
{\bf x}_{k}=\bm{\Phi}\bm{\Lambda}^k{\bf q}_0. 
\end{equation}
Or equivalently, it is written as 
\begin{equation}
{\bf x}_k=\sum_{i=1}^{m}\mu_i^k\bm{\varphi}_i q_{0,i}, \label{eq:rev:2}
\end{equation}
where $q_{0,i}$ denotes the $i$th component of ${\bf q}_0$. 
Now consider a continuous-time dynamical system behind the discrete dynamical system defined by \eqref{eq:th:-1}, i.e., 
\begin{equation}
\dot{\bf x}=\bm{\mathcal{A}}{\bf x}, \label{eq:rev:4}
\end{equation}
where $\bm{\mathcal{A}}$ is the continuous counterpart of ${\bf A}$. 
The solution of \eqref{eq:rev:2} to a given initial condition ${\bf x}(0)$ can be obtained as (see, e.g.,~\cite{Chen2009}), 
\begin{equation}
{\bf x}(t)=\sum_{i=1}^{m}{\rm e}^{s_i t}\bm{\phi}_ip_{0,i}, \label{eq:rev:3}
\end{equation}
where $s_i$ and $\bm{\phi}_i$ denote the $i$th eigenvalue and eigenvector of $\bm{\mathcal{A}}$, $p_{0,i}$ is the $i$th modal coordinate of ${\bf p}_0=\bm{\Psi}^*{\bf x}(0)$ where $\bm{\Psi}=[\bm{\phi}_1,\dots,\bm{\phi}_m]$. 
By comparing Eqs.~(\ref{eq:rev:2}) and (\ref{eq:rev:3}) at $t=t_k=k\Delta t$ with the assumption that $\bm{\phi}_i\approx\bm{\varphi}_i$ and $p_{0,i}\approx q_{0,i}$, one see that 
\begin{equation}
{\rm e}^{s_ik\Delta t}=\mu_i^k. 
\end{equation}
This yields
\begin{equation}
s_i=\log(\mu_i)/\Delta t\label{eq:th:0}. 
\end{equation}
As indicated by \eqref{eq:rev:3}, the eigenvalue $s_i$ contains the temporal characteristics of $\bm{\phi}_i$, which is typically represented by using natural frequency and damping ratio, as follows, 
\begin{equation}
s_i=-\zeta_i\omega_i\pm{\rm j}\omega_i\sqrt{1-\zeta_i^2}, 
\end{equation}
where $\omega_i$ denotes the natural frequency and $\zeta_i$ denotes the corresponding modal damping ratio. 
Therefore, undamped natural frequencies $f_i$, and modal damping ratios $\zeta_i$ can be extracted as follows. 
\begin{align}
f_i&=|s_i|/2\pi\label{eq:th:1},\\
\zeta_i&=-{\rm Re}(s_i)/|s_i|\label{eq:th:2}. 
\end{align}
%
}}
It is known that the DMD yields an approximate eigendecomposition of the best-fit linear operator relating two data matrices~\cite{Tu2013}, which in this case is the matrix ${\bf A}$. 
%
Note also that it is known that the DMD produces the approximations of {\it Koopman operator eigenfunctions}~\cite{WilliamsEtAl2015}. 
Furthermore, normal modes of the vibrating system and the Koopman operator eigenfunctions extracted from time series of the system are known to be identical if the system is linear~\cite{CirilloEtAl2016}. 
Therefore, for vibrating systems, since the eigenvalues and eigenvectors of ${\bf A}$ {\color{black}{approximate}} the eigenvalues and eigenvectors of the underlying dynamical system{\color{black}{, or $\bm{\mathcal{A}}$ of \eqref{eq:rev:4}}}, the DMD modes and their eigenvalues should correspond to their normal modes and corresponding eigenvalues, which contain undamped natural frequencies and corresponding damping ratios. 

%
Since the process of interest in this paper is linear, 
the DMD modes to be obtained have the nature of standing waves. 
It is known that the DMD procedure defined above 
fails to capture the behavior of standing waves. For such cases, the snapshot matrices {\bf X} and {\bf Y} need to be redefined with time-shifted snapshots appended, as follows~\cite{TuEtAl2014}, 
\begin{align}
{\bf X}&\triangleq
\begin{pmatrix}
{\bf x}_1,\dots,{\bf x}_{n\hspace{1em}}\\
{\bf x}_2,\dots,{\bf x}_{n+1}
\end{pmatrix},\label{eq:th:sw1}\\
{\bf Y}&\triangleq
\begin{pmatrix}
{\bf x}_2,\dots,{\bf x}_{n+1}\\
{\bf x}_3,\dots,{\bf x}_{n+2}\label{eq:th:sw2}
\end{pmatrix}. 
\end{align}
The DMD computations are then applied to Eqs.~(\ref{eq:th:sw1}) and (\ref{eq:th:sw2}). 
\subsection{Relationship with the Ibrahim Time Domain method}
In this subsection, the Ibrahim Time Domain method (ITD) is briefly reviewed and its relationship with the DMD is discussed. 
The ITD is a time-domain modal parameter extraction method, where a unique set of modal parameters are obtained from a set of free vibration measurements in a single measurement~\cite{Ewins2009}. The ITD assumes the free response of the viscously-damped MDOF system of the form:
\begin{equation}
{\bf M}\ddot{\bf u}(t)+{\bf C}\dot{\bf u}(t)+{\bf K}{\bf u}(t)={\bf 0}, \label{eq:th:4}
\end{equation}
where ${\bf M}$, ${\bf C}$, and ${\bf K}$ denote mass, damping and stiffness matrices, ${\bf u}(t)$ is a vector of displacement coordinates, ${\bf M},{\bf C},{\bf K}\in\mathbb{R}^{m\times m}$, and ${\bf u}(t)\in\mathbb{R}^m$, $m$ is the number of DOF of the system. 
At a given time instant $t=t_j$ where $t_j\in[t_1,t_{n}]$, ${\bf u}(t_j)$ can be represented by a linear combination of the eigenvectors of the system \eqref{eq:th:4}, i.e., 
\begin{equation}
{\bf u}(t_j)=
\sum_{i=1}^{2n}\bm{\psi}_i{\rm e}^{s_it_j}\label{eq:th:5}, 
\end{equation}
where $s_i$ denotes the $i$th eigenvalue, $\bm{\psi}_i$ denotes the corresponding eigenvector. 
Note that $s_i$ and $\bm{\psi}_i$ contain $m$ of complex conjugate pairs. 
Defining matrices ${\bf U}_{2n}=[{\bf u}(t_1),\dots,{\bf u}(t_{2n})]$, and $\bm{\Phi}=[\bm{\phi}_1,\dots,\bm{\phi}_{2n}]$, 
the following relationship holds. 
\begin{equation}
{\bf U}_{2n}=\bm{\Phi}\bm{\Lambda}\label{eq:th:55}, 
\end{equation}
where 
\begin{equation}
\bm{\Lambda}=
\begin{pmatrix}
{\rm e}^{s_{1}t_1} &\cdots &{\rm e}^{s_1t_{2n}}\\
\vdots         & & \vdots\\
{\rm e}^{s_{2n}t_{1}} &\cdots &{\rm e}^{s_{2n}t_{2n}}
\end{pmatrix}. 
\end{equation}
Note that ${\bf U}\in\mathbb{R}^{n\times2n}$, $\bm{\Phi}\in\mathbb{C}^{n\times2n}$ and $\bm{\Lambda}\in\mathbb{C}^{2n\times2n}$. 
Next consider the following expansion at the time instant $t_{j+1}$ where $t_{j+1}=t_{j}+\Delta t$ for $j=1,\dots,2n$, 
\begin{equation}
{\bf u}(t_{j+1})=\sum_{i=1}^{2n}\tilde{\bm{\phi}}_i{\rm e}^{s_it_j}, 
\end{equation}
where $\tilde{\bm{\phi}}_i=\bm{\phi}_i{\rm e}^{s_i\Delta t}$. Since this holds for $t_{j+1}$ for $i=1,\dots,2n$, the following relationship is obtained, i.e., 
\begin{equation}
{\bf U}_{2n+1}=\tilde{\bm{\Phi}}\bm{\Lambda}\label{eq:th:6}, 
\end{equation}
where ${\bf U}_{2n+1}=[{\bf u}(t_2),\dots,{\bf u}(t_{2n+1})]$, $\tilde{\bm{\Phi}}=[\tilde{\bm{\phi}}_1,\dots,\tilde{\bm{\phi}}_{2n}]$. 
Finally, the third set of equations are considered for a time instant $t_{j+2}$ where $t_{j+2}=t_{j}+2\Delta t$ for $j=1,\dots,2n$, 
\begin{equation}
{\bf u}(t_{j+2})=\sum_{i=1}^{2n}\hat{\bm{\phi}}_i{\rm e}^{s_it_j}, 
\end{equation}
where $\hat{\bm{\phi}}_i=\bm{\phi}_i{\rm e}^{s_i2\Delta t}$.
\begin{equation}
{\bf U}_{2n+2}=\hat{\bm{\Phi}}\bm{\Lambda},\label{eq:th:7}
\end{equation}
where ${\bf U}_{2n+2}=[{\bf u}(t_3),\dots,{\bf u}(t_{2n+2})]$ and $\hat{\bm{\Phi}}=[\hat{\bm{\phi}}_1,\dots,\hat{\bm{\phi}}_{2n}]$. 
From Eqs~(\ref{eq:th:5}), (\ref{eq:th:6}), and (\ref{eq:th:7}), we obtain
\begin{align}
{\bf X}&=\bm{\Psi}\bm{\Lambda}, \\
\tilde{\bf X}&=\tilde{\bm{\Psi}}\bm{\Lambda}, 
\end{align}
where 
\begin{align*}
{\bf X}=\begin{bmatrix}
{\bf U}_{2n\hspace{1em}}\\{\bf U}_{2n+1}
\end{bmatrix},
\tilde{\bf X}=\begin{bmatrix}
{\bf U}_{2n+1}\\{\bf U}_{2n+2}
\end{bmatrix},
\bm{\Psi}=\begin{bmatrix}
\bm{\Phi}\\
\tilde{\bm{\Phi}}
\end{bmatrix},
\tilde{\bm{\Psi}}=\begin{bmatrix}
\tilde{\bm{\Phi}}\\
\hat{\bm{\Phi}}
\end{bmatrix}.
\end{align*}
This leads to 
\begin{align}
\tilde{\bm{\Psi}}&=\tilde{\bf X}\bm{\Lambda}^{-1}\nonumber\\ 
&=\tilde{\bf X}{\bf X}^{-1}{\bm{\Psi}}. \label{eq:th:8}
\end{align}
Denoting the $i$th columns of $\bm{\Psi}$ and $\tilde{\bm{\Psi}}$ respectively as $\bm{\psi}_i$ and $\tilde{\bm{\psi}}_i$, the following relationship holds:
\begin{equation}
\tilde{\bm{\psi}}_i=\bm{\psi}_i{\rm e}^{s_i\Delta t}\label{eq:th:9}. 
\end{equation}
Therefore, substituting \eqref{eq:th:9} into \eqref{eq:th:8} yields the following eigenvalue problem:
\begin{equation}
\left(\tilde{\bf X}{\bf X}^{-1}\right)\bm{\psi}_i=\mu_i\bm{\psi}_i,\quad i=1,\dots,{2n},  \label{eq:th:10}
\end{equation}
where $\mu_i={\rm e}^{s_i\Delta t}$, and $\mu_i\in\mathbb{C}$. If the number of time instants in the snapshot ${\bf U}$ is greater than the system size $2n$, then the inversion of ${\bf X}$ in \eqref{eq:th:8} is not possible because ${\bf X}$ is rectangular. Therefore, the inverse matrix is replaced with the Moore-Penrose pseudo inverse matrix, as follows, 
\begin{equation}
\tilde{\bm{\Psi}}=\tilde{\bf X}{\bf X}^{\dagger}\bm{\Psi}, 
\end{equation}
where ${\bf X}^{\dagger}={\bf X}^{\rm T}\left({\bf X}{\bf X}^{\rm T}\right)^{-1}$. The eigenvalue problem \eqref{eq:th:10} then becomes as follows:
\begin{equation}
\left(\tilde{\bf X}{\bf X}^{\dagger}\right)\bm{\psi}_i=\mu_i\bm{\psi}_i,\quad i=1,\dots,2n. 
\label{eq:th:11}
\end{equation}
By definition, the first $n$ entries in $\bm{\psi}_i$ contain the system's $i$th eigenvector. Therefore, by computing $\mu_i$ and $\bm{\psi}_i$ from the snapshot matrices ${\bf X}$ and $\tilde{\bf X}$, the modal parameters of the system can be obtained. 

In the context of DMD, the matrix ${\bf A}$ defined in \eqref{eq:th:-1+1} is equivalent to $\tilde{\bf X}{\bf X}^{\dagger}$ for the ITD method. Therefore, the eigenvalues and eigenvectors obtained by solving \eqref{eq:th:11} are equivalent to the DMD eigenvalues and eigenvectors. 
However, in DMD, since the pseudo inverse of ${\bf A}$ is obtained via the SVD, the singular values with small magnitudes can be neglected along with their left and right singular vectors. 
This means that the pseudo-inversion involved in DMD is more numerically stable than that used in ITD method. 
Furthermore, the DMD modes can be ranked and ordered by the magnitude of their singular values. 
%
{\color{black}{
\subsection{Least squares complex frequency domain method}
The LSCF method is one of the most widely used frequency-domain modal parameter extraction methods. The mathematical formulation of the LSCF method is briefly reviewed. 
Denoting the FRF of $o$th output point of the object as $H_o(\omega)$, the LCSF method assumes that it can be represented as~\cite{VanDerAuweraerEtAl2001,PeetersEtAl2004}, 
\begin{equation}
H_o(\omega)=B_o(\omega)/A(\omega), \label{eq1:lscf}
\end{equation}
and 
\begin{align}
B_o(\omega)&=\sum^{N_p}_{p=0}b_{o,p}{\rm e}^{-{\rm j}p\omega T_s},\\
A(\omega)&=\sum^{N_p}_{p=0}a_{p}{\rm e}^{-{\rm j}p\omega T_s}, 
\end{align}
where $N_p$ denotes the maximum polynomial order to be considered, $b_{o,p}$ and $a_p$ represent unknown coefficients to be determined, $T_s$ is the sampling period. The unknown coefficients are estimated by a weighted linear least squares method. Grouping the unknown coefficients as $\bm{\theta}=[b_{o,0},\cdots,b_{o,N_p}, a_{0},\cdots,a_{N_p}]$, the weighted linear least squares method is defined as 
\begin{equation}
\minimize_{\bm{\theta}} A(\omega_k)\varepsilon_o(\omega_k,\bm{\theta}),\quad\mbox{for }k=1,\dots,N_f
\end{equation}
where $N_f$ denotes the number of frequency lines and $\varepsilon_o$ is defined as the weighted differece between the measured FRF $\hat{H}_o$ and its approximation represented by \eqref{eq1:lscf}, i.e., 
\begin{equation}
\varepsilon_o(\omega_k,\bm{\theta})=w(\omega_k)\left[
H_o(\omega_f,\bm{\theta})-\hat{H}_o(\omega_f)
\right],\quad\mbox{for }k=1,\dots,N_f, 
\end{equation}
where $w(\omega_k)$ is the frequency-dependent weight that can be a function of coherence function, for instance. 
This equation is solved for a given value of $N_p$, and the eigenvalues of the characteristic polynomial are identified as the poles of the system, which include the information of natural frequencies and corresponding damping ratios of the structure. 
The drawback of this procedure is that it generates mathematical poles that do not represent physical vibration modes. 
To circumvent this problem, this process is repeated with increasing $N_p$, and the {\it stabilization diagram} can be drawn, where the poles that appear regardless of $N_p$ are identified as stable poles, i.e., physical vibration modes. }}
\section{Numerical examples}\label{section:numerical_example}
This section provides numerical examples for the application of the DMD to identify modal properties of mechanical structures. 
\subsection{Single DOF system}
\begin{figure}[tb]
\centering
\includegraphics[scale=.7]{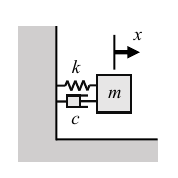}
\caption{Single DOF oscillator model}\label{fig:sdof:model}
\end{figure}
Consider a single DOF damped oscillator shown in \figref{fig:sdof:model}, where $m=1$kg, $c=6.2832$N$\cdot$s/m, and $k=9.8696\times10^4$N/m. For this particular example, the chosen parameter values yield natural frequency $f_n=\omega_n/(2\pi)=50$Hz, $\zeta=0.01$, and damped natural frequency $f_d=49.9997$Hz. 
Since this is an underdamped 1DOF oscillator, the time response of the mass for a given set of initial conditions $x(0)=x_0$ and $\dot{x}(0)=v_0$ can be written as~\cite{Meirovitch2001}, 
\begin{equation}
x(t)={\rm e}^{-\zeta\omega_nt}\left(
x_0\cos(\omega_d t) + \frac{v_0+\zeta\omega_nx_0}{\omega_d}\sin(\omega_d t)
\right), \label{eq1}
\end{equation}
where $\omega_d=\omega_n\sqrt{1-\zeta^2}$ is the damped natural frequency, $\omega_n=\sqrt{k/m}$ is the undamped natural frequency and $\zeta=c/2\sqrt{mk}$ is the damping ratio.  
\begin{figure}[tb]
\centering
\subfigure[Time history]{\includegraphics{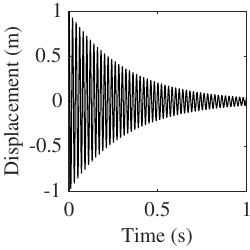}}
\subfigure[DMD eigenvalues]{\includegraphics{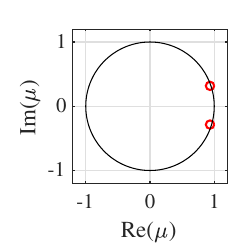}}
\caption{Time history of the mass and extracted DMD eigenvalues}\label{fig:sdof:1}
\end{figure}
The displacement of the mass for the initial conditions of $x_0=1$ and $v_0=0$ has been computed for $0\leqslant t \leqslant 1$ using \eqref{eq1}, and sampled at $n=1024$ discrete time instants. The time response of the displacement is shown in \figref{fig:sdof:1}(a). The displacement shows damped oscillation as expected. 
The DMD is then applied to the obtained snapshots based on the procedure described in \secref{section:theory}. 
In \figref{fig:sdof:1}(b), the obtained DMD eigenvalues $\mu$ are shown in a complex plane with a unit circle.  
As can be seen, the DMD eigenvalues appear as a complex conjugate pair, where 
\begin{equation}
\mu=0.9503\pm0.3014{\rm j}.
\end{equation}
They are slightly inside the unit circle due to the damping, i.e., $|\mu|=0.9969$. 
Application of Eqs.~(\ref{eq:th:1}) and (\ref{eq:th:2}) yields $f_n$ and $\zeta$. 
The obtained results are shown in \tabref{tab:sdof:1}. 
As can be seen, $f_n$ and $\zeta$ extracted by the DMD algorithm agree very well with the exact values. 
\begin{table}[tb]
\centering
\caption{Extracted natural frequency and damping ratio by DMD}\label{tab:sdof:1}
\begin{tabular}{c|cc}\hline
& $f_n$ & $\zeta$\\ \hline
Exact & $5.00\times10^1$Hz & $1.00\times10^{-2}$ \\ 
DMD & $4.99\times10^1$Hz & $1.00\times10^{-2}$ \\ 
Relative error & $1.13\times10^{-13}\%$ & $2.41\times10^{-12}\%$\\ \hline
\end{tabular}
\end{table}
\subsection{Six degrees of freedom system}\label{subsec:sixdof}
\begin{figure}[tb]
\centering
\includegraphics{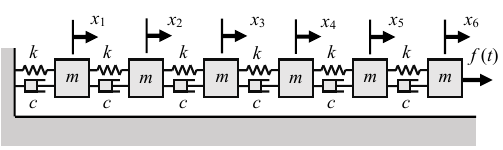}
\caption{Six degree of freedom system}\label{fig:mdof:1}
\end{figure}
Next, as a representative example of multi-DOF mechanical systems, a six DOF mass-spring-damper system shown in \figref{fig:mdof:1} is considered. The equations of motion of the system can be written as a matrix-vector representation of the form: 
\begin{equation}
{\bf M}\ddot{\bf x}(t)+{\bf C}\dot{\bf x}(t) + {\bf K}{\bf x}(t)={\bf f}(t), \label{eq:mdof:1}
\end{equation}
where ${\bf x}=[x_1,x_2,x_3,x_4,x_5,x_6]^{\rm T}$, ${\bf M}$, {\bf C}$, {\bf K}$ are mass, damping and stiffness matrices, and they are defined as
\begin{align}
{\bf M}&=
\begin{pmatrix}
m & 0 & 0 & 0 & 0 &0\\
 0 & m & 0 & 0 & 0 & 0\\
 0 & 0 & m & 0 & 0 & 0\\
 0 & 0 & 0 & m & 0 & 0\\
 0 & 0 & 0 & 0 & m & 0\\
 0 & 0 & 0 & 0 & 0 & m
\end{pmatrix},\\
{\bf K}&=
\begin{pmatrix}
2k & -k & 0 & 0 & 0 &0\\
 -k & 2k & -k & 0 & 0 & 0\\
 0 & -k & 2k & -k & 0 & 0\\
 0 & 0 & -k & 2k & -k & 0\\
 0 & 0 & 0 & -k & 2k & -k\\
 0 & 0 & 0 & 0 & -k & k
 \end{pmatrix},
\end{align}
and ${\bf C}=(c/k){\bf K}$. ${\bf f}(t)$ is the external forcing vector where ${\bf f}(t)=[0,0,0,0,0,f(t)]^{\rm T}$. 
Following the standard mode superposition procedure, given the mass-normalized modal matrix $\bm{\Phi}=[\bm{\phi}_1,\dots,\bm{\phi}_6]$, the equation of motion \eqref{eq:mdof:1} can be rewritten as, 
\begin{equation}
\ddot{z}_i(t)+2\zeta_i\omega_i\dot{z}_i(t)+\omega_i^2z_i(t)=g_i(t), 
\end{equation}
where $z_i(t)$ is the $i$-th modal coordinate that satisfies 
\begin{equation}
{\bf x}(t)=\bm{\Phi}{\bf z}(t), \label{eq:mdof:1-2}
\end{equation}
where ${\bf z}(t)=[z_1(t),\dots,z_6(t)]^{\rm T}$, $\omega_i$ is the $i$th (undamped) natural frequency, $\zeta_i$ is the $i$th damping ratio, $g_i(t)$ is the $i$th modal forcing defined as $g_i(t)=\bm{\phi}_i^{\rm T}{\bf f}(t)$. 
In this numerical example, consider a case where a unit step forcing is applied to the 6th DOF, i.e., $f(t)=1$ for $t\geqslant 0$. 
Also, $m=1$, $c=0.01$, and $k=1$ are used for the generation of the numerical example. This yields the natural frequencies and corresponding modal damping ratios shown in \tabref{tab:mdof:1}. 
\begin{table}[tb]
\centering
\caption{Modal damping ratios and undamped natural frequencies of the six degrees of freedom system}\label{tab:mdof:1}
\begin{tabular}{c|cc}\hline
$i$ & $f_{n,i}$ (Hz) & $\zeta_i$ \\ \hline
1& $0.0384$ & $1.2054\times10^{-3}$ \\
2& $0.1129$ & $3.5460\times10^{-3}$ \\
3& $0.1808$ & $5.6806\times10^{-3}$ \\
4& $0.2383$ & $7.4851\times10^{-3}$ \\
5& $0.2818$ & $8.8546\times10^{-3}$ \\
6& $0.3091$ & $9.7094\times10^{-3}$ \\ \hline
\end{tabular}
\end{table}
As can be seen in \tabref{tab:mdof:1}, all modal coordinates are underdamped. Therefore, analytic solutions of $z_i(t)$ can easily be found, as follows. 
\begin{equation}
z_i(t)=
\frac{g_i}{\omega_i^2}\left[
1-{\rm e}^{-\zeta_i\omega_{i}t}\left(
\cos\omega_{di}t+\frac{\zeta_i}{\sqrt{1-\zeta_i^2}}\sin\omega_{di}t
\right)
\right]\label{eq:mdof:2}, 
\end{equation}
where $\omega_{di}=\omega_i\sqrt{1-\zeta_i^2}$. The physical DOFs can be recovered by \eqref{eq:mdof:1-2}. 
\begin{figure}[bt]
\subfigure[Time response of modal coordinates (${\bf z}(t)$)]{\includegraphics{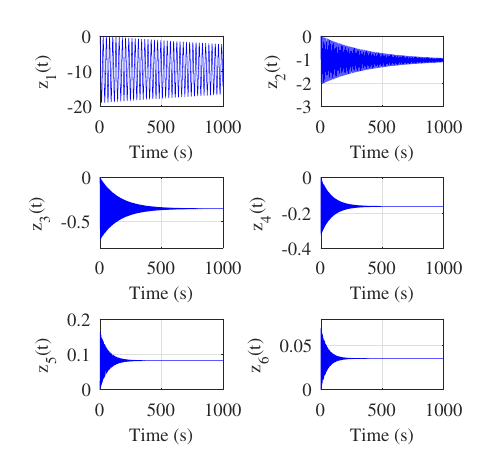}}
\subfigure[Time response of physical DOF (${\bf x}(t)$)]{\includegraphics{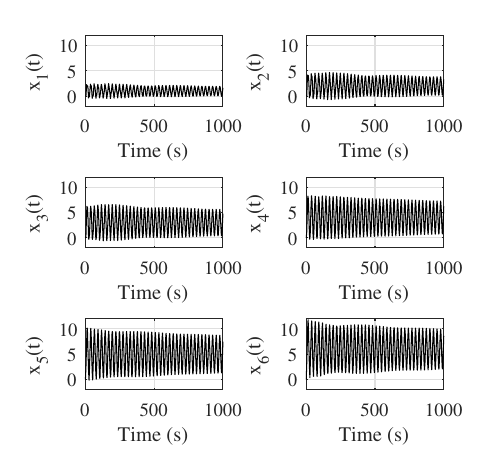}}
\caption{Time response of modal and physical DOFs}\label{fig:mdof:2}
\end{figure}

The time responses of modal coordinates were computed based on \eqref{eq:mdof:2} for 
$0\leqslant t \leqslant 1000[{\rm s}]$, and 
are shown in \figref{fig:mdof:2}(a). 
As can be seen, as the mode index increases, the frequency of the modal coordinates increase and they tend to converge quicker, as the damping ratio increases as mode index increases as well, as shown in ~\tabref{tab:mdof:1}. 
The recovered time histories of ${\bf x}$ are shown in \figref{fig:mdof:2}(b). Obviously the dominant frequency component for all DOFs is $f_{n,1}$. 
Since the step force is applied to the sixth DOF, the center and the amplitude of oscillation of the sixth DOF become the largest among the DOFs. These values decrease as the DOF index decreases, or as the DOF approaches to the fixed end. 
\begin{figure}[bt]
	\subfigure[{\color{black}{DMD and displacement spectra}}: $\times$, DMD, {\color{blue}{---}}, Displacement]{\includegraphics{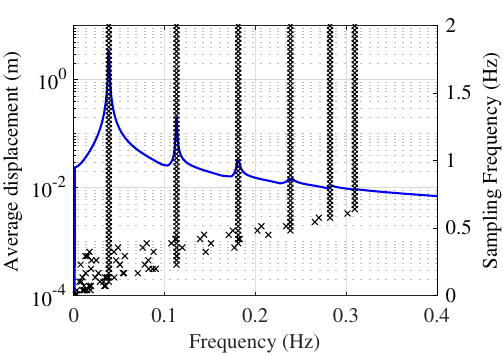}}
	\subfigure[DMD damping ratio: $\times$, DMD, {\color{red}{---}}, exact]{\includegraphics{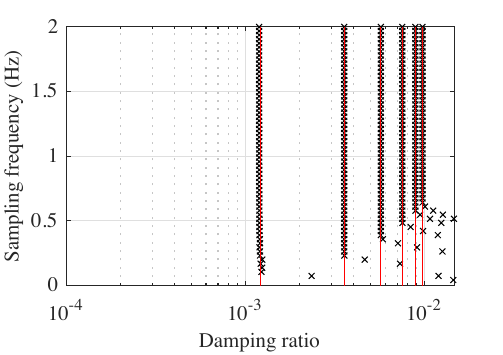}}
	\caption{Pseudo-stability diagrams computed by DMD for the six DOF system}\label{fig:mdof:3}
\end{figure}

The modal parameters have then been extracted by the DMD procedure described in \secref{section:theory}. The obtained DMD spectra are shown in \figref{fig:mdof:3}. 
To examine the effects of the sampling frequency on the results, the results are shown with the increasing value of the sampling frequency. 
This representation of the natural frequencies and damping ratios with the corresponding sampling frequency is hereinafter referred to as pseudo-stability diagram in this paper, because it resembles the representation of the natural frequencies in the existing modal parameter extraction methods, such as LSCF. 
Figure~\ref{fig:mdof:3}a shows the undamped natural frequencies extracted by DMD, combined with {\color{black}{displacement spectrum.}}
%
As expected, correct natural frequency is captured after the sampling frequency exceeds approximately twice the natural frequency. 
Namely, a natural frequency that is larger than the Nyquist frequency of the sampling frequency can be captured by the DMD. 
Figure~\ref{fig:mdof:3}b shows the corresponding damping ratios for increasing values of sampling frequency. 
The vertical lines in the graph show the exact values. 
As can be seen, accurate damping ratios tend to be captured after the sampling frequency exceeds twice the corresponding natural frequencies. 
Overall, one can see that both natural frequency and damping ratios can be extracted accurately with the DMD. 
{\color{black}{It is noted that even though many modes appear to be captured for each frequency, they indeed represent the same mode. The criterion employed here is that if a pole consistently appears regardless of the sampling frequency, then it represents a physical mode, which is inspired by the concept of the stability diagram used in the LSCF method. Therefore, the pole with any of the sampling frequencies can be selected as long as it lies within a predefined threshold, which can be set based on the average of the poles computed with all sampling frequencies considered, for instance. }}
\begin{figure}[bt]
	\centering
	\includegraphics{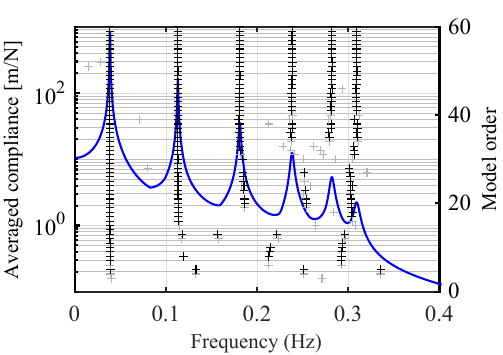}
	\caption{Stability diagram computed by LSCF method: $+$, stable poles with negative real parts, 
	${\color{gray}{+}}$: unstable pole,
	$$\textcolor{blue}{---}, averaged compliance}\label{fig:mdof:4}
\end{figure}

To compare the DMD with the existing modal parameter extraction method, the LSCF method~\cite{VanDerAuweraerEtAl2001} has been applied to the FRF. 
The order of polynomial has been increased up to 60, and the stability diagram is obtained as shown in \figref{fig:mdof:4}. {\color{black}{The results are shown for polynomial orders from 1 to 60. }}
{\color {black}{In \figref{fig:mdof:4}, it is noted that the black cross indicates a stable pole that appears as a pole even with increasing polynomial order within a pre-defined threshold, which in this case is chosen to be 1\%. 
The grey crosses indicate poles that do not stay within the threshold, which means that they are not stable and considered to be artificial poles. }}
First, the accuracy of the natural frequencies and damping ratios are discussed. As can be seen in \figref{fig:mdof:4}, the LSCF method can also capture the natural frequencies accurately. 
\begin{table}[tb]
\centering
\caption{Percentage errors in natural frequencies and damping ratios: DMD is computed with $f_s=2$Hz. LSCF is computed with model order of 60. }\label{tab:mdof:2}
\begin{tabular}{c|cc|cc}\hline
          & \multicolumn{2}{c|}{Frequency}& \multicolumn{2}{c}{Damping ratio}\\ \hline
$i$     & DMD & LSCF & DMD & LSCF \\ \hline
1 & 8.079$\times10^{-5}$ & 1.671$\times10^{-3}$ & 1.500 & 2.530 \\ 
2 & 1.771$\times10^{-5}$ & 6.193$\times10^{-3}$ & 1.216$\times10^{-2}$ & 6.965$\times10^{-1}$ \\ 
3 & 5.530$\times10^{-6}$ & 1.456$\times10^{-2}$ & 3.128$\times10^{-3}$ & 4.376$\times10^{-1}$ \\ 
4 & 1.259$\times10^{-5}$ & 2.783$\times10^{-2}$ & 3.295$\times10^{-3}$ & 8.399$\times10^{-1}$ \\ 
5 & 2.128$\times10^{-5}$ & 2.937$\times10^{-2}$ & 3.384$\times10^{-3}$ & 7.145$\times10^{-1}$ \\ 
6 & 9.706$\times10^{-6}$ & 2.110$\times10^{-2}$ & 3.965$\times10^{-3}$ & 7.876$\times10^{-2}$ \\ \hline
\end{tabular}
\end{table}
The results of both DMD and LSCF methods are tabulated in \tabref{tab:mdof:2} for both natural frequencies and damping ratios. 
Both methods can capture both natural frequencies and damping {\color {black}{ratios}} well for all modes. 
The DMD tends to capture the frequency much better than the LSCF method. For the damping ratios, even though the difference is only about an order of magnitude, DMD appears to capture the values better than the LSCF. 
{\color {black}{In experimental modal analysis, it is commonly preferred to use time-domain method for systems with a low damping ($\zeta<1\%$), while frequency-domain methods, such as LSCF, are preferred for highly damped systems. This is due mostly to the fact that FRF data from low damped systems contains most of its energy in a very limited number of frequency lines, whereas the time-domain data has a large number of time samples~\cite{Verboven2002}. }}
These results indicate that the DMD method can extract both natural frequencies and damping ratios as good as the existing frequency-domain modal parameter extraction method. 

{\color {black}{It is noted that the LSCF method typically requires more data, because the computation of FRF requires the information of the input, which in this case is the force applied to the system. On the other hand, DMD does not require the information of the input.}}
\begin{figure}[tb]
\centering
\subfigure[Mode 1]{\includegraphics[scale=.9]{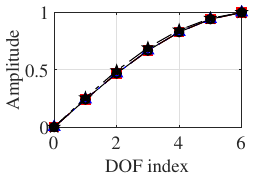}}
\subfigure[Mode 2]{\includegraphics[scale=.9]{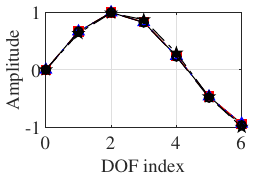}}
\subfigure[Mode 3]{\includegraphics[scale=.9]{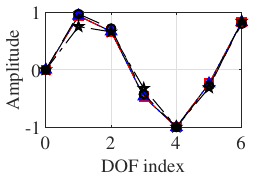}}
\subfigure[Mode 4]{\includegraphics[scale=.9]{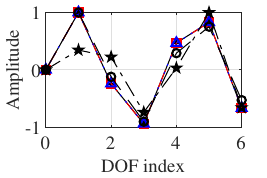}}
\subfigure[Mode 5]{\includegraphics[scale=.9]{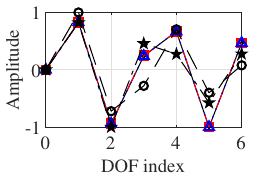}}
\subfigure[Mode 6]{\includegraphics[scale=.9]{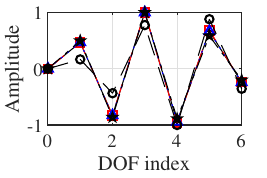}}
\caption{Mode shapes. Analytic, ---+---, LSCF, {\color{red}{---\footnotesize{$\square$}---}}, DMD (\#1\&\#2), {\color{blue}{$-\!-\!\!\!\triangle\!\!\!-\!-$}}, POD\#1, $-\!-\!\!\circ\!\!-\!-$, POD\#2, --- - --$\star$-- - ---}\label{fig:mdof:5}
\end{figure}
\begin{table}[tb]
\centering
\caption{Comparison of MAC values (in \%) computed with the analytic solutions. }\label{tab:mdof:3}
\begin{tabular}{c|cccc}\hline
Mode index & LSCF & DMD & POD\#1 & POD\#2\\ \hline
1 & 100.00 & 100.00 & 100.00 &  99.97 \\ 
2 & 100.00 & 100.00 &  99.97 &  99.79 \\ 
3 & 100.00 & 100.00 &  99.78 &  98.38 \\ 
4 & 100.00 & 100.00 &  98.39 &  72.75 \\ 
5 & 100.00 & 100.00 &  73.97 &  88.03 \\ 
6 &  99.99 & 100.00 &  89.27 &  99.63 \\ \hline
\end{tabular}
\end{table}

Next, the accuracy of the mode shapes extracted by the DMD method is discussed. The DMD mode shapes $\bm{\varphi}$ have been computed by \eqref{eq:th:32}, and shown in \figref{fig:mdof:5}. 
The DMD modes appear as complex conjugate pairs. Therefore, only real parts of the DMD modes are shown in \figref{fig:mdof:5} (because the imaginary part of a DMD mode is proportional to its real part). 
In \figref{fig:mdof:5}, analytic solutions of the mode shapes, the ones extracted by the LSCF (hereinafter LSCF modes), and the POD modes are also shown for comparison. 
There are two POD modes that resemble each linear vibration mode. Therefore, the two POD modes are labeled as POD\#1 and \#2, where the singular value corresponding to POD\#1 is greater than that corresponding to POD\#2. 
As can be seen, DMD, LSCF, and both POD modes are almost identical for the first two modes. 
However, from Mode 3, the differences between the POD modes and the other modes become apparent. 
As the mode order increases, the differences appear to increase. 
%
%
It is known that the POD modes converge to the linear normal modes as the amount of data increases~\cite{FeenyKappagantu1998}. Therefore, these discrepancies may disappear if more data points were taken into account. 
On the other hand, both LSCF and DMD methods show excellent agreement with the analytic solutions. 
In addition to the mode shapes themselves, the MAC values (see, e.g., ~\cite{Ewins2009}) between DMD, LSCF, POD modes and the analytic solutions have been computed and shown in \tabref{tab:mdof:3}. 
As shown in the table, both DMD and LSCF produced accurate mode shapes for all six modes. 
The MAC values decrease as the mode index increases for the POD modes, as it was also visible in \figref{fig:mdof:5}. 
\subsection{Effects of measurement errors}\label{sec:error}
\begin{figure}[tb]
\subfigure[{\color{black}{DMD and displacement spectra}} for $\sigma=1.0\times10^{-5}$: $\times$, DMD, {\color{blue}{---}}, Without noise, {\color{red}{-$\cdot$-}}, With noise]{\includegraphics{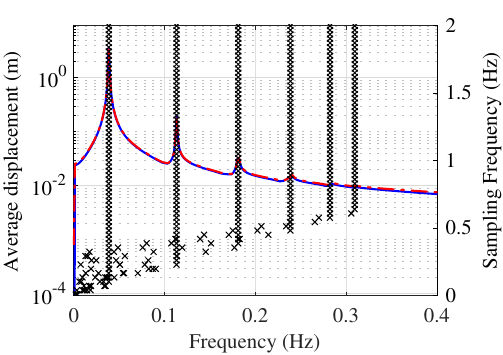}}
\subfigure[DMD damping ratio for $\sigma=1.0\times10^{-5}$: $\times$, DMD, {\color{red}{---}}, exact]{\includegraphics{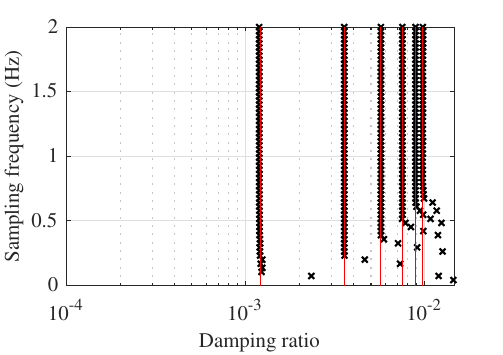}}
\subfigure[{\color{black}{DMD and displacement spectra}} for $\sigma=1.0\times10^{-2}$: $\times$, DMD, {\color{blue}{---}}, Without noise, {\color{red}{-$\cdot$-}}, With noise]{\includegraphics{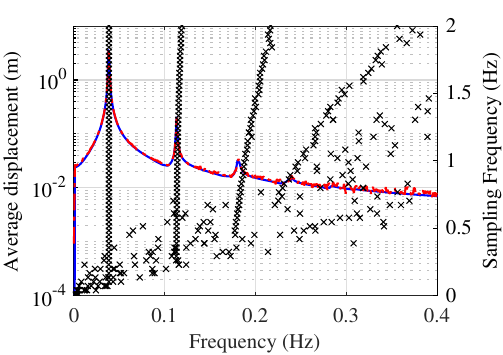}}
\subfigure[DMD damping ratio for $\sigma=1.0\times10^{-2}$: $\times$, DMD, {\color{red}{---}}, exact]{\includegraphics{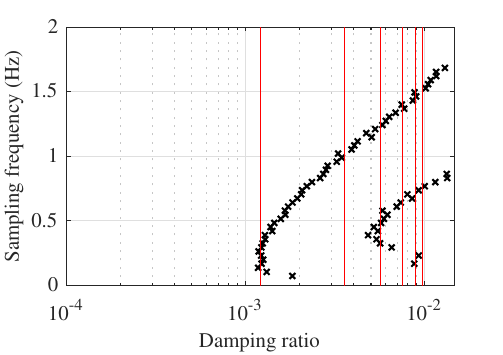}}
\caption{Pseudo-stability diagrams computed by DMD for the six DOF system with measurement errors}\label{fig:error:1}
\end{figure}
Next, effects of measurement errors in the time response on the results of DMD is examined, by using the numerical example used in \secref{subsec:sixdof}. 
The time response of the system subjected to a unit step forcing on the 6th DOF is again computed based on Eqs.~(\ref{eq:mdof:1}) and (\ref{eq:mdof:2}). This time, in order to study the effects of measurement errors on the DMD results, artificial measurement errors are injected into the time histories of the DOFs, as follows:
\begin{equation}
\tilde{x}_i(t_j)=\left[1+
\sigma\delta_i(t_j)
\right]
x_i(t_j),
\quad\mbox{for }i = 1,\dots,6, \quad j=1,\dots,n_t, 
\end{equation}
where $n_t$ is the number of time instants in the snapshots, and $\delta_i(t_j)$ are a set of uniformly distributed random numbers with unit standard deviation and zero mean, and $\sigma$ is the pre-defined standard deviation. 
The DMD was then applied to the data and modal parameters are extracted for $\sigma=1.0\times10^{-5}$ and $\sigma=1.0\times10^{-2}$, with the increasing sampling frequency up to 2Hz, as it was done in \secref{subsec:sixdof}. 
The results of the DMD are shown in \figref{fig:error:1}. 
As can be seen in \figref{fig:error:1}(a) and (b), with relatively small measurement errors with $\sigma=1.0\times10^{-5}$, DMD can extract both natural frequencies and damping ratios with high accuracy. 
On the other hand, with relatively large measurement errors with $\sigma=1.0\times10^{-2}$, DMD cannot extract both natural frequencies and damping ratios. In particular, third to sixth natural frequencies cannot be well captured for the increasing sampling frequency. 
It is worth noting that the damping ratios cannot be well captured by the DMD, with the presence of this level of measurement errors. 
These results align with the discussion by Duke {\it et al.}~\cite{DukeETAl2012}, where it is shown that the grow or decay rate predicted by DMD is highly sensitive to the signal to noise ratios.
This high sensitivity to the measurement errors can be an issue to be overcome especially when the DMD is applied to noisy experimental data and damping ratio is of critical interest. 

{\color{black}{In order to further examine the sensitivity of the DMD eigenvalues on the measurement noise, transmission of measurement errors on the DMD eigenvalues is considered. Let us take the first order perturbations on $s_i$ and $\mu_i$ as $s_i'=s_i+\delta s_i$ and $\mu_i'=\mu_i+\delta \mu_i$, where $\delta s_i$ and $\delta\mu_i$ represent perturbations on the eigenvalues caused by measurement errors in ${\bf x}$, or ${\bf A}$. By substituting $s_i'$ and $\mu'_i$ into \eqref{eq:th:0}, we obtain, 
\begin{equation}
s_i+\delta s_i = f_s\log\left(\mu_i+\delta \mu_i\right),\label{eq:sensitivity}
\end{equation}
where $f_s=1/\Delta t$ is the sampling frequency. 
By taking advantage of \eqref{eq:th:0}, we obtain 
\begin{equation}
\delta s_i = f_s\log\left(
1+\delta \mu_i/\mu_i
\right).
\end{equation}
This means that if $|\delta\mu_i/\mu_i|\ll1$, $|\delta s_i|\approx f_s|\delta\mu_i/\mu_i|$. Since $|\delta\mu_i/\mu_i|$ is finite for given measurement error, this means that $|\delta s_i|$ is approximately proportional to $f_s$. This argument aligns with the observations for \figref{fig:error:1}, where one sees that increasing the sampling frequency decreases the prediction accuracy. }}

Next, the effects of measurement errors on the mode shapes are discussed. The MAC values were computed between the DMD modes and the analytically obtained mode shapes, and shown in \tabref{tab:error:1}. 
As can be seen in \tabref{tab:error:1}, for $\sigma=1.0\times10^{-5}$, the MAC values for both DMD and POD modes are as accurate as the one without measurement errors, which are shown in \tabref{tab:mdof:3}. 
On the other hand, for $\sigma=1.0\times10^{-2}$, even though the first DMD mode still keeps relatively good MAC value of 99.93\%, MAC values for the higher modes are low, which means the mode shapes extracted by DMD do not resemble the analytic solutions. 

Again, this high sensitivity of DMD needs to be improved especially when noisy experimental data are used. 
{\color {black}{It is known that the application of total least squares method by introducing augmented snapshot matrix can be effective in treating measurement errors in DMD~\cite{HematiEtAl2017}, however, it is beyond the scope of this paper. }}
\begin{table}[tb]
\centering
\caption{Comparison of MAC values (in \%) for the mode shapes obtained by DMD with noisy data. }\label{tab:error:1}
\begin{tabular}{c|ccc|ccc}\hline
& \multicolumn{3}{c}{$\sigma=1.0\times10^{-5}$}&\multicolumn{3}{|c}{$\sigma=1.0\times10^{-2}$}\\ \hline
Mode index & DMD & POD\#1 & POD\#2 &DMD & POD\#1 & POD\#2\\ \hline
1 & 100.00 & 100.00 &  99.97  &  99.93 & 100.00 & 100.00 \\ 
2 & 100.00 &  99.97 &  99.79   &  59.43 &  99.99 &  99.44 \\ 
3 & 100.00 &  99.78 &  98.38   &  32.78 &  99.83 &  86.90 \\ 
4 & 100.00 &  98.39 &  72.75   &  71.64 &  90.34 &   9.70 \\ 
5 & 100.00 &  73.97 &  88.03   &  25.34 &  48.85 &   5.85 \\ 
6 & 100.00 &  89.27 &  99.63   &  68.36 &  62.84 &  29.95 \\ \hline
\end{tabular}
\end{table}
%
%
%
\section{Experimental example: impulse response of a cantilever beam}\label{section:experimental_example}
\begin{figure}[bt]
\centering
\subfigure[Test specimen]{\includegraphics[width=8.5cm]{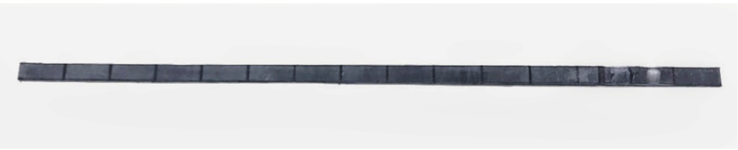}}
\subfigure[Experimental setup]{\includegraphics[width=8.5cm]{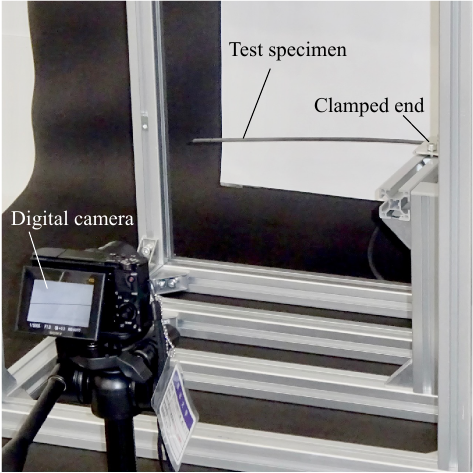}}
\caption{Test specimen and experimental setup}\label{fig:experimental_setup}
\end{figure}
In this section, results of the application of DMD to experimental data are discussed. Of particular interest is the extraction of modal parameters from image-based large-scale data obtained from high-speed camera. 
In this study, vibration response of a slender cantilevered beam subjected to impulse excitation is used for the application of DMD. 
{\color{black}{In a relevant development of DMD, the application of DMD to the response for impulsive loadings has been reported in \cite{KouEtAl2017}, where the dynamic modes for subcritical flows are successfully extracted.}} 
Figure \ref{fig:experimental_setup} shows the test specimen and the experimental setup used for this study. The beam is made of polypropylene and its dimension is 300$\times$8.0$\times$1.7mm with rectangular cross sectional area (Fig.\ref{fig:experimental_setup}(a)). 
The beam was attached to a fixture and an impulsive force was applied near the fixed end by an impulse hammer (GK-2110, Ono Sokki, Japan). 
The resulting vibration of the beam was recorded as a sequence of images by a high-speed digital camera (RX100M4, Sony, Japan) and saved as a single movie file. 
The frame rate of the movie was 
480 frame per second (fps), which should be enough for capturing the vibration modes up to 240 Hz. 
The color information contained at each pixel in the images in the movie was binarized so that the edge of the beam can be detected. 
The time history of the displacement of the test specimen was then extracted from the images along the axial direction in a range of 15.64mm to 204.2mm measured from the free end (to the left in \figref{fig:experimental_setup}(b)), at equally-spaced 1,570 discrete sampling points.  
The extraction of the images from the movie file, binarization of color information, and conversion of pixel to physical length were all conducted by Matlab\textsuperscript{\textregistered}. 
\begin{figure}[bt]
\centering
{\color{black}{
\subfigure[$t=1.41$s]{\includegraphics[width=8cm]{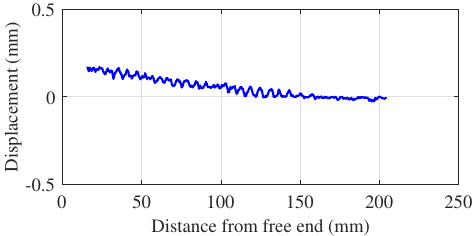}}
\subfigure[$t=1.47$s]{\includegraphics[width=8cm]{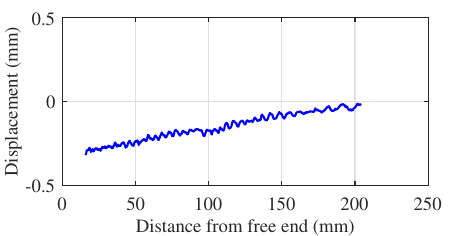}}
\subfigure[$t=1.52$s]{\includegraphics[width=8cm]{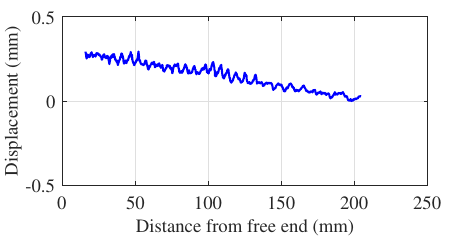}}
\subfigure[$t=1.57$s]{\includegraphics[width=8cm]{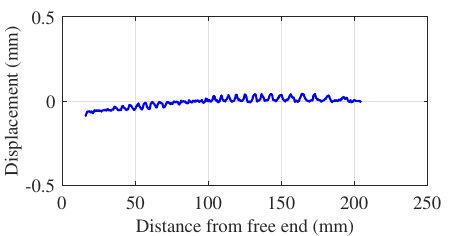}}
}}
\caption{Measured vibration response of the beam {\color{black}{at representative time instants}}
}\label{fig:impact_frame}
\end{figure}

In \figref{fig:impact_frame}, measured 
{\color{black}{displacement }}
of the beam after binarization at representative time instances are shown. Note that the time instant when the impact of the impulse hammer occurred is defined to be $t=0$. 
As can be seen, the impact of the impulse hammer induced the elastic vibration of the cantilevered beam. 
As expected, multiple vibration modes appears to have been excited by the impulsive forcing. 
{\color {black} Also, as can be seen, the response shapes contain some measurement noise caused by the binalization of the image. }
\begin{figure}[tb]
\centering
\includegraphics[scale=.9]{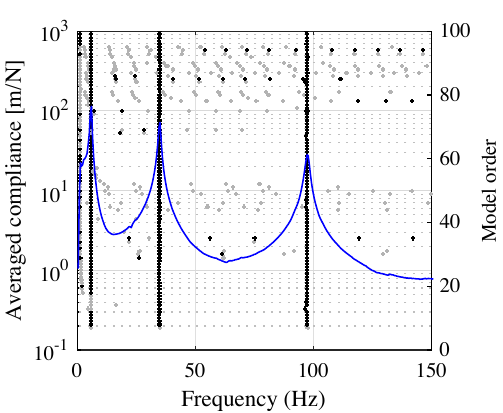}
\caption{Stability diagram obtained from LSCF method with 1,570 measurement points. $+$: stable pole with negative real part, $\filledcirc$: unstable pole with positive real part, \textcolor{blue}{---}: averaged compliance.}\label{exp:fig:lscf}
\end{figure}
The time histories of the displacement were extracted for the 1,570 points on the beam, and DMD has been applied to the time histories. As it was done in \secref{subsec:sixdof}, DMD computation has been repeated for the increasing sampling frequency for 200Hz$\leqslant f_s\leqslant$ 450Hz. 
Furthermore, in order to compare the results obtained by the DMD with the ones obtained from LSCF method, FRFs and coherence functions have also been computed by using the spectra of the input force and the obtained displacement field.  
{\color{black}{For obtaining the FRFs, five measurements were conducted and the average of the FRFs in frequency domain was used for the LSCF method.}} 

First, results obtained by LSCF method are discussed. The result is shown in \figref{exp:fig:lscf}, where the averaged compliance and identified natural frequencies are shown with increasing polynomial order. It is noted that it is known that the image-based FRFs require higher-order polynomials in the LSCF method~\cite{JavhEtAl2018}. 
Thus, the model order was increased up to 100. 
As can be seen, by using LSCF method, three modes were identified below 150Hz with clearly stable poles.  {\color {black}{It is noted that the LSCF method has detected an unstable pole near 0Hz. Even though the pole was detected as a stable pole with the highest polynomial order, the pole was unstable with polynomials with lower orders. This pole is detected because of a small “edge” in the FRF data near 0Hz.}}
\begin{figure}[tb]
\centering
\subfigure[DMD spectra combined with averaged Fourier spectrum of displacement. $\times$, DMD, {\color{blue}{---}}, Fourier spectrum of displacement]{\includegraphics[scale=.8]{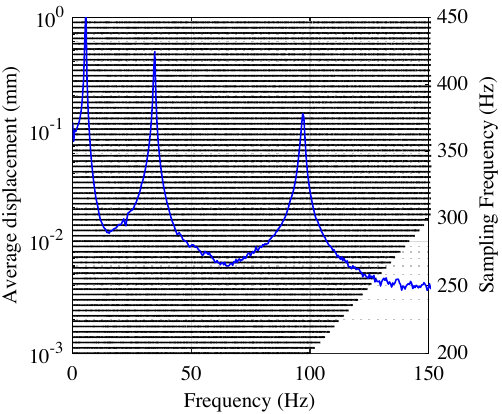}}
\subfigure[DMD damping ratio]{\includegraphics[scale=.8]{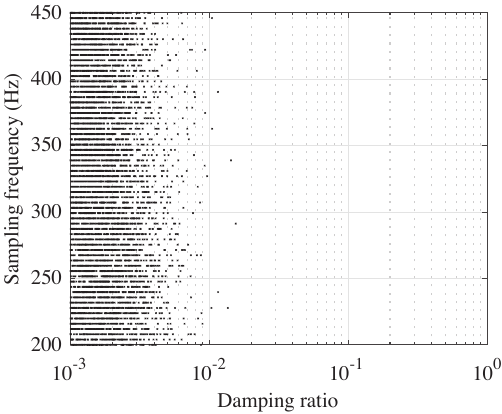}}
\caption{Pseudo-stability diagrams computed by DMD without singular value rejection}\label{exp:fig:dmd1}
\end{figure}
\begin{figure}
\centering
\subfigure[DMD spectra combined with averaged autospectrum of displacement. $\times$, DMD, {\color{blue}{---}}, Fourier spectrum of displacement]{\includegraphics[scale=.8]{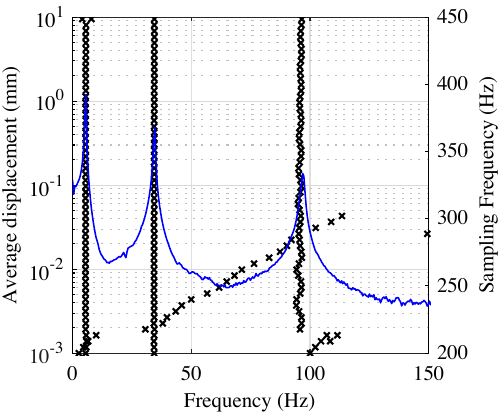}}
\subfigure[DMD damping ratio]{\includegraphics[scale=.8]{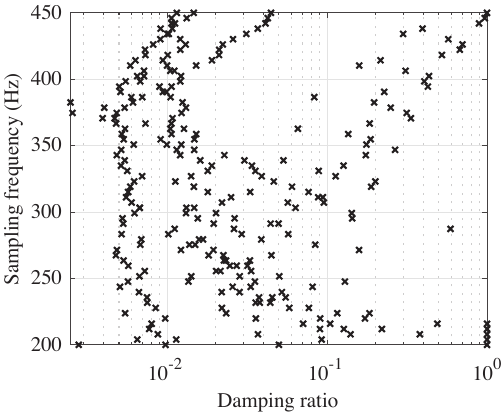}}
\caption{Pseudo-stability diagrams computed by DMD with singular value rejection}\label{exp:fig:dmd2}
\end{figure}

Next, results obtained by DMD are discussed. 
In \figref{exp:fig:dmd1}, pseudo-stability diagrams of natural frequencies and damping ratios are shown with increasing sampling frequency. 
As can be seen, both natural frequencies and damping ratios predicted by the DMD produce almost continuous spectra. With these diagrams, it is not feasible to identify the real modal parameters of the system. 
This is because the dataset used for DMD contains measurement noise with many measurement points. 
Therefore, DMD produced many eigenvalues, resulting in almost continuous spectra. 
This requires a strategy to filter out insignificant eigenvalues produced from DMD. 
{\color{black}{There are many methods to date to select appropriate dynamic modes, such as the ones based on the evolution of the dynamic modes in sampling space~\cite{KouZhang2017}}}. 
One of the most effective approaches to circumvent this problem is the use of singular value rejection (e.g., \cite{ThiteThompson2003}), where pseudo-inverse matrix is constructed with significant singular values and singular vectors, by ignoring the singular values with small magnitude and corresponding left and right singular vectors. 
%
%
In \figref{exp:fig:dmd2}, pseudo-stability diagrams of natural frequencies and damping ratios are presented, where the singular value rejection has been applied with the DMD. One can see that insignificant poles have all been removed from the diagrams, resulting in clear stability diagram especially for the natural frequencies. Again, damping ratios corresponding to the identified natural frequencies do not produce stable loci in the diagram, which was also observed in \secref{sec:error}. 

{\color{black}{It should be mentioned that the results shown in \figref{exp:fig:dmd1} can be considered as the results of ITD method where singular value rejection is not applied. This case study shows that the DMD can produce numerically more stable results than the traditional ITD method. }}

\begin{table}[tb]
\centering
\caption{Comparison between LSCF and DMD in terms of natural frequencies and damping ratios}\label{exp:tab:f_and_zeta}
\begin{tabular}{c|ccc|ccc}\hline
& \multicolumn{3}{c|}{Frequency [Hz]}& \multicolumn{3}{c}{Damping ratio}\\ \hline
Mode index & DMD & LSCF & Error [\%] & DMD & LSCF &Error [\%]\\ \hline
1 & 5.91 & 5.61 & 5.08 & 0.0460 & 0.06290 & 36.14\\
2 & 34.90 & 34.41 & 1.39 & 0.01429 & 0.00731 & 48.86\\
3 & 97.13 & 96.35 & 0.81 & 0.01181 & 0.01100 & 6.84\\ \hline
\end{tabular}
\end{table}
First three modes were identified by both methods, and the obtained natural frequencies and damping ratios are tabulated in \tabref{exp:tab:f_and_zeta}. 
Note that the results of DMD with sampling frequency of 450Hz, and the results of LSCF with the polynomial order of 100 have been used for the comparison. 
As can be seen, natural frequencies obtained from DMD agree well with the ones obtained from LSCF method, with up to 6\% error. 
On the other hand, damping ratios obtained from DMD do not agree well with those obtained from LSCF. 

It is noted that in order to confirm the validity of the results obtained by DMD, numerical modal analysis has also been conducted by creating a three dimensional beam model discretized by finite element method (FEM), by using ANSYS\textsuperscript{\textregistered}. The dimension of the model is $250\times8.0\times1.7$mm, where Young's modulus is 1.6GPa, Poisson's ratio is 0.4, and density is 900kg/m$^3$. The beam model was discretized with solid elements with the edge length being set as 0.6mm. 
One of the faces was rigidly fixed to simulate the clamped end. 
The obtained first three natural frequencies were 5.88Hz, 36.8Hz, and 103Hz, which are in good agreement with those obtained by both DMD and LSCF method. 
Note that the natural frequencies are slightly larger than those experimentally obtained by DMD and LSCF methods, mostly because of the boundary condition employed, because rigidly fixing a face typically yields stiffer system than the real system with clamped end. 
\begin{figure}[tb]
\centering
\subfigure[Mode 1]{\includegraphics[scale=.9]{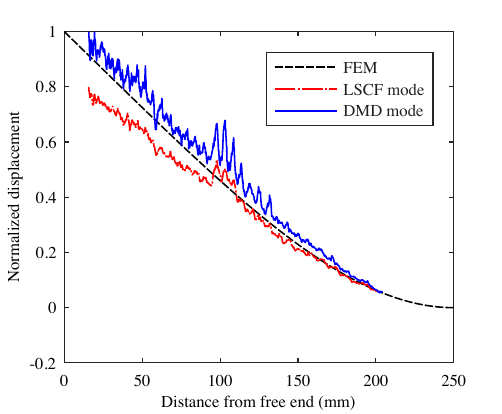}}
\subfigure[Mode 2]{\includegraphics[scale=.9]{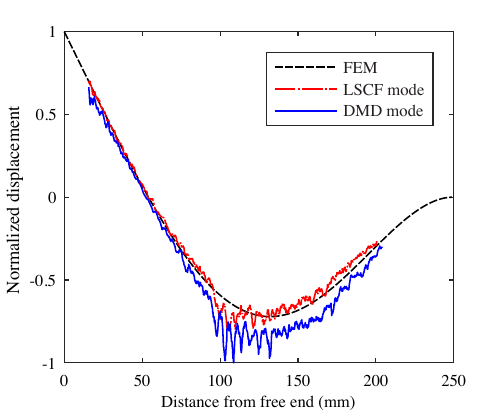}}
\subfigure[Mode 3]{\includegraphics[scale=.9]{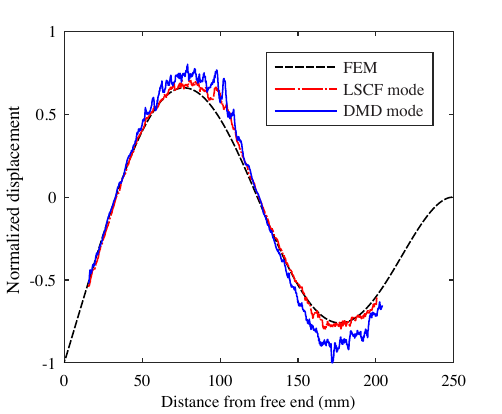}}
\caption{Mode shapes extracted from measured data by LSCF, DMD, and computed by FEM}\label{fig:experimental_example:1}
\end{figure}

\begin{table}[tb]
\centering
\caption{MAC (in \%) between the modes computed by DMD and FEM}\label{table:experimental_example:1}
\begin{tabular}{c|ccc}\hline
\diagbox{DMD}{FEM}
&1 & 2 & 3\\ \hline
1  &   99.56  &  8.365 &   6.130 \\
2 &   10.48  & 98.85  &  2.609\\
3  &    6.602   & 2.336  & 99.39 \\ \hline
\end{tabular}
\end{table}
\begin{table}[tb]
\centering
\caption{MAC (in \%) between the modes computed by LSCF and FEM}\label{table:experimental_example:2}
\begin{tabular}{c|ccc}\hline
\diagbox{LSCF}{FEM}
&1 & 2 & 3\\ \hline
1  &     99.57  &  8.960  &  6.053 \\ 
2 &       6.559  & 99.37  &  1.750\\
3  &      6.524  &  1.443 &  99.52\\ \hline
\end{tabular}
\end{table}

Next, mode shapes obtained by DMD, LSCF, and FEM are presented and their differences are discussed. In \figref{fig:experimental_example:1}, mode shapes obtained by DMD, LSCF, and FEM are presented for the first three modes. 
As can be seen in the figure, both DMD and LSCF produced mode shapes that agree well with the ones obtained by FEM. 
{\color {black}{One can observe that there are small noises in the mode shapes for both LSCF and DMD. This is because of the noise in the snapshot matrix, as it was shown in \figref{fig:impact_frame}. The noise level in the LSCF modes appears to be slightly lower than that in the DMD modes. This is attributed to the fact that the averaging technique employed in the FRF measurement, which is typically done for removing measurement noise in EMA. On the other hand, the snapshot matrix used for the computation of DMD is from a single measurement. The issue of removing such spatial noise is an issue to be resolved in the future work.}} 
To quantitatively evaluate the accuracy of the mode shapes predicted by the method, MAC values have been computed. 
MAC values between the modes computed by DMD and FEM are shown in \tabref{table:experimental_example:1}. 
As can be seen in the diagonal elements in the table, the mode shapes computed by DMD and FEM are all in good agreement with each other, with MAC values being over 98\%. 
Furthermore, as the off-diagonal terms show, the orthogonality of the modes is also confirmed. 
MAC values between the modes computed by LSCF and FEM are also shown in \tabref{table:experimental_example:2}. 
In the table, similar trends can be found with slightly better MAC values on the diagonal elements. 
Overall, however, the results are consistent with the ones obtained from DMD. 
\section{Conclusions}\label{section:conclusion}
In this paper, a data-driven experimental modal analysis based on DMD was examined. 
%
DMD shares the fundamental concept with a well-known Ibrahim time domain method. However, in DMD, the DMD eigenvalues are computed based on the transformation matrix that is projected onto the POD modes, with insignificant POD modes being rejected. 
%
First, DMD was applied to analytically-obtained time histories of discrete mass-spring-damper systems. 
With the absence of measurement errors, DMD can capture the modal parameters very accurately. 
Second, DMD was applied to the same time history with measurement noise. 
When the standard deviation of measurement errors is small, DMD can still capture the modal parameters accurately. 
However, when the standard deviation of the measurement errors is relatively large, DMD fails to capture the modal parameters. 
%
Furthermore, DMD was applied to time histories of displacement of cantilevered beam acquired by high-speed digital camera. The obtained modal parameters were compared with the ones obtained by LSCF method. The natural frequencies and corresponding mode shapes obtained by DMD with the application of singular value rejection agreed well with the ones obtained by LSCF method. The damping ratios did not agree well with those obtained by LSCF method. 

In summary, DMD has a potential as a method for EMA especially when a large amount of data is involved. However, our study indicates that the estimation of damping ratios accurately still remains a challenge when applying DMD to noisy experimental data. 

\end{document}